# A Logic-Based Mixed-Integer Nonlinear Programming Model to Solve Non-Convex and Non-Smooth Economic Dispatch Problems— An Accuracy Analysis

M. Pourakbari-Kasmaei[1*], M. Fotuhi-Firuzabad[2], and J. R. Sanches Mantovani[3]

[1] Department of Electrical Engineering, Aalto University, Maarintie 8, 02150 Espoo, Finland.
[2] Department of Electrical Engineering, Sharif University of Technology, Center of Excellence in Power System Management & Control (CEPSMC), Tehran, Iran.
[3] Department of Electrical Engineering, State University of Sao Paulo, Ilha Solteira, Brazil.
* Mahdi.Pourakbari@aalto.fi

*Abstract:* This paper presents a solver-friendly logic-based mixed-integer nonlinear programming model (LB-MINLP) to solve economic dispatch (ED) problems considering disjoint operating zones and valve-point effects. A simultaneous consideration of transmission losses and logical constraints in ED problems causes difficulties either in the linearization procedure, or in handling via heuristic-based approaches, and this may result in outcome violation. The non-smooth terms can make the situation even worse. On the other hand, non-convex nonlinear models with logical constraints are not solvable using the existing nonlinear commercial solvers. In order to explain and remedy these shortcomings, we proposed a novel recasting strategy to overcome the hurdle of solving such complicated problems with the aid of the existing nonlinear solvers. The proposed model can facilitate the pre-solving and probing techniques of the commercial solvers by recasting the logical constraints into the mixed-integer terms of the objective function. It consequently results in a higher accuracy of the model and a better computational efficiency. The acquired results demonstrated that the LB-MINLP model, compared to the existing (heuristic-based and solver-based) models in the literature, can easily handle the non-smooth and nonlinear terms and achieve an optimal solution much faster and without any outcome violation.

## 1. Introduction

Electricity is one of the most important infrastructures of any developed/developing country in the world. However, as the nations' economic and social welfare improves, the global energy consumption, and particularly the electricity demand, increases. On the other hand, the energy cost has sharply increased in recent years. This reveals the main role of the economic analyses in power system from tie-line planning [1], generation and transmission expansion planning [2,3], distribution network expansion planning [4–7], toward power system operation studies such as unit commitment [8–10], maintenance scheduling [11], optimal power flow [12,13], and economic dispatch (ED) [14,15], etc. In this paper we discuss the key role of ED tool in an optimal distribution of electricity demands among all the committed units [16]. It is an appropriate tool to determine the optimal output of all the committed units at a minimum cost, while satisfying the hourly demands and the operational constraints [17]. A simple ED model can consider a continuous and convex cost function. However, unlike the traditional ED models, in practice, the models are non-convex and non-smooth because of the prohibited operating zones (POZs), the valve-point effects, and the multi-fuel options [18].

In the literature, various mathematical models and optimization approaches have been used to solve the ED problems, aiming at obtaining the optimal solution with less computational time. These optimization approaches can be classified for modeling the language- and the heuristic-based optimization techniques [19]. The heuristic-based approaches are useful for solving complicated problems that are not solvable via the existing commercial solvers. Among the proposed heuristic-based approaches in the literature, the dimensional steepest decline (DSD) method is one of the most successful heuristic-based approaches for solving complicated ED problems [20]. Although this approach succeeded to consider the non-smooth terms, for the sake of simplicity, the linear and constant parts of the transmission loss in the Kron's loss formula were neglected, which resulted in an inappropriate model. However, for all the heuristic-based approaches, due to the requirement of several trials and adjustments for each system and under each condition, the precision of the solution is still not guaranteed. Such a lack of guarantee shows the importance of investigating the language-based models, which are solvable by commercial solvers, namely the solver-based models, which due to the rapidity and precision of finding the optimal solution, are currently the most attracted approaches.

The first solver-based model that addressed the POZs was presented in [21]. It was an interesting mixed-integer quadratic programming (MIQP) model for finding the optimal solution of the ED problems with disjoint operating zones. However, the main drawback of this model was using too many auxiliary decision-making variables that were assigned to achieve this model. For each unit with *k*-1 POZs, *k* continuous variables for the feasible zones and k binary variables to jump over the POZs were assigned. The second drawback of this MIQP model was the incapability of considering the transmission losses. To address this problem, which is related to the use of too many auxiliary variables, the big-M based MIQP (M-MIQP) was proposed in [22]. This model, by using only $k$-1 artificial binary variables via a binary coding approach, achieved the same optimal solution with a higher computational efficiency. However, besides the lack of considering the transmission losses, it conveyed two new disadvantages as follows: 1) Using the big-M method to consider the POZs; in fact, finding a proper M is a trial-and-error task, and a very large M yields serious numerical difficulties in the computer; 2) Complex binary encoding. Later, in [23], the unambiguous distance-based MIQP model (UDB-MIQP) addressed the two drawbacks of the M-MIQP model by assigning only $k$-1 binary variables to the distances among the operating zones, which consequently obtained the same optimal results with much higher computational efficiency. However, the lack of considering the transmission losses was still a shortage. It is worth mentioning that none of the above MIQP-based models are capable of considering the non-smooth terms.

Although the aforementioned methods are very fast and have brought a new vision to solving the ED problems, the lack of considering the transmission losses and the non-smooth terms are their undeniable shortcomings. It should be expressed that transmission losses may vary from 5 to 15 percent of the total load. Therefore, representation of the transmission losses in the economic dispatch calculations is indispensable. On the other hand, in practice, non-smoothness and non-convexity are the inseparable characteristics of power systems. These deficiencies have motivated different researchers to introduce new approaches that are capable of addressing such practical constraints. Recently, a mixed-integer quadratically constrained quadratic programming (MIQCQP) model was suggested to solve the shortcoming of transmission losses [24]. The MIQCQP model exposed violation in the power balance, which might be due to the transmission loss miscalculation. Moreover, this model, while having the same two drawbacks as in [22], conveyed two more problems as follows: 1) Using a complex bi-level procedure for solving a normally single-level problem; and 2) Using different solvers for the upper and lower levels. Therefore, the MIQCQP model failed to consider the transmission losses in an accurate way by showing violation in some cases.

All in all, from the recently performed works and the published results, it can be deduced that the transmission loss, with or without considering the non-smooth terms, is still an important issue in ED problems to be addressed. This paper aims at addressing the aforementioned drawbacks by proposing a logic-based mixed-integer nonlinear programming (LB-MINLP) model. The proposed model by making a perfect tradeoff between the model accuracy and the computational efficiency tries to reach the best outcome without losing the concept.

The contributions of the current paper are as follows:

1) Filling the existing gap in the literature related to the nonlinear and non-smooth

characteristics of power systems, by considering the transmission loss calculated using the Kron's loss formula, and proposing a logic-based MINLP model.

2) Proposing a new recasting approach that brings facilities for commercial nonlinear solvers. Via this recasting, the logical constraints of the MINLP model are reformulated into mixed-integer terms and are inserted into an objective function. This recasting facilitates the pre-solving and probing techniques of commercial solvers and results in a higher computational efficiency and a more accurate output.

3) Performing the accuracy and the power balance violation checks for the proposed model and the existing approaches in the literature. The transmission losses of the works in the literature were re-calculated with the definition of the Kron's loss formula in [25] in order to demonstrate either the accuracy or the misinterpretation of the Kron's loss formula.

In order to verify the proposed LB-MINLP model and reveal its strengths and weaknesses, the commonly used 6- and 15-unit test systems are considered in detail. Then, in order to show the effectiveness of the proposed LB-MINLP in solving the practical large-scale power systems, the practical Korean 140- and 10,000-unit power systems are investigated.

The remainder of this paper is organized as follows. Section II presents the proposed LB-MINLP. The case studies and results are demonstrated in Section III. Section IV provides the concluding remarks and the prospects for future works.

## 2. The LB-MINLP Framework for Non-convex and Non-smooth Economic Dispatch Problems

In this section, at first, the mathematical model of a practical ED problem is presented and then, the proposed LB-MINLP is explained in more detail.

### 2.1. The mathematical formulation of a practical ED problem

The main objective of an ED problem is to minimize the total generation cost, as follows.

$$\min \sum_{i=1}^{n_g} F_i(P_i) \tag{1}$$

where $i$, $n_g$, and $P_i$ denote the index of each generator, the total number of generators, and the output power of generator $i$, respectively. $F_i(\cdot)$ signifies the generation cost of unit $i$, which is mostly approximated by a quadratic function, as shown in (2):

$$F_i(P_{g_i}) = a_i(P_i)^2 + b_i P_i + c_i \tag{2}$$

where $a_i$, $b_i$, and $c_i$ are the quadratic, linear, and constant cost coefficients of unit $i$, respectively.

On the other hand, it should be stated that in practice, multiple valves will result in the ripples. In turbines, these valves are designed to control the steam for separating the groups of nozzles in the first stage of the turbine and to ultimately achieve a higher average efficiency, over a wide range of loads by successively admitting steam to the groups of nozzles and meet the demand fluctuations [26]. Therefore, it is significantly essential to consider the valve-points in ED problems. The valve-point effect can be modeled by adding a sinusoid term in the cost function (2), as described below [14].

$$F_i(P_{g_i}) = a_i(P_i)^2 + b_i P_i + c_i + \left| e_i \times \sin(f_i \times (\underline{P_i} - P_i)) \right| \tag{3}$$

where $e_i$ and $f_i$ are the cost coefficients, corresponding to the valve-point effect. In addition, $\underline{P_i}$ is the minimum generation capacity of unit $i$.

In ED problems, several constraints must be satisfied, as explicated in the following:

*1) Power balance constraint*

In a power system, the total power demand and the transmission losses must be met.

$$\sum_{i=1}^{n_g} P_i = P_D + P_L \tag{4}$$

where $P_D$ is the power demand and $P_L$ is the transmission loss, calculated via the B-coefficient method [25], as shown in (5).

$$P_L = \sum_{i=1}^{n_g} \sum_{j=1}^{n_g} P_i B_{ij} P_j + \sum_{i=1}^{n_g} B_{0,i} P_i + B_{00} \tag{5}$$

where $B_{ij}$, $B_{0i}$, and $B_{00}$ are, respectively, the $ij^{th}$ element of the loss coefficient square matrix, the $i^{th}$ element of the loss coefficient vector, and the loss coefficient constant. It is worth mentioning that in order to calculate the system loss, in Equation (5), all the output powers must be considered per unit [27].

*2) Active power output limit*

Each generator must generate between its lower and upper output capacity limits.

$$\underline{P_i} \leq P_i \leq \overline{P_i} \tag{6}$$

where $\overline{P_i}$ is the maximum generation capacity of unit $i$.

*3) Spinning reserve*

The spinning reserve is the online reserve capacity of all the committed units, which is used during the emergency operating conditions, outages, contingencies, and unforeseen load swings, to maintain the system frequency stability.

$$\sum_{i=1}^{n_g} S_i \geq S_R \tag{7}$$

$$S_i = \min\{(\overline{P}_i - P_i), \overline{S}_i\} \tag{8}$$

where $S_i$ and $\overline{S}_i$ are the available and the maximum spinning reserves of unit $i$, respectively. $S_R$ stands for the required spinning reserve of the system.

*4) Ramp-rate limit*

Each generator has its ramp-up (9) and ramp-down limits (10), which must be satisfied during the increase or the decrease of the output power, respectively.

$$P_i - P_i^0 \leq UR_i \tag{9}$$

$$P_i^0 - P_i \leq DR_i \tag{10}$$

Therefore, by considering Equations (6), (9), and (10), the output of each generator must satisfy the following constraint:

$$\max(\underline{P}_i, P_i^0 - DR_i) \leq P_i \leq \min(\overline{P}_i, P_i^0 + UR_i) \tag{11}$$

where $P_i^0$ is the previously scheduled output power of unit $i$. In addition, $DR_i$ and $UR_i$ are the ramp-down and ramp-up rates of generator $i$, respectively.

*5) Disjoint operating zones*

In practice, due to some faults in the machine or its accessories, or due to the vibrations in the shaft bearing (caused by the steam valve), the generating unit might have prohibited the operating zones, which can be formulated as (12) [28].

$$\begin{cases} \underline{P}_i = \underline{P}_{i1} \leq P_{i1} \leq \overline{P}_{i1}, \text{ or} \\ \underline{P}_{ik} \leq P_{ik} \leq \overline{P}_{ik}, \quad 2 \leq k \leq (z_i - 1), \text{ or} \\ \underline{P}_{iz_i} \leq P_{iz_i} \leq \overline{P}_{iz_i} = \overline{P}_i \end{cases} \tag{12}$$

where $k$ denotes the index of the operating zones. Moreover, $\overline{P}_{ik}$ and $\underline{P}_{ik}$ stand for the upper and lower bounds of the operating zone $k$ of unit $i$. $z_i$ is the total number of the (feasible) operating zones of unit $i$.

### 2.2. The logic-based mixed-integer nonlinear programming framework

The main idea of the proposed framework comes from the fact that each generator can only generate in one of its operating zones, and each operating zone has its own sub-cost function. Therefore, this framework is a sub-cost function selecting-based model. For further illustration, the process of selection, for a single unit $i$ with two POZs, was considered in detail. It should be stated that for the sake of simplicity, in the rest of this paper, the sub-function is used instead of the sub-cost function.

In Figure 1, the cost function of the generating unit, considering two POZs (the regions with red hachures), is depicted. The solid and dashed curves are, respectively, related to the quadratic cost function (2) and the cost function with considering valve-points (3). This unit has three operating zones and three corresponding sub-cost functions ($F_{i1}$, $F_{i2}$, and $F_{i3}$).

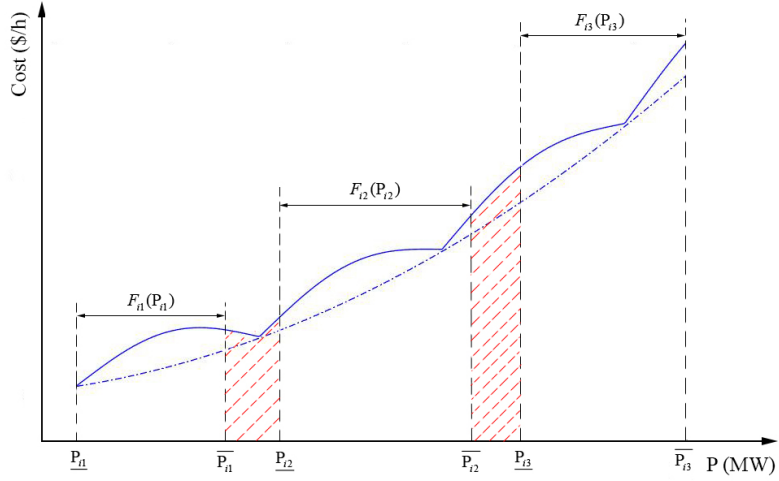

Figure 1. Valve point effect of unit $i$ with three feasible zones (two POZs)

$$\begin{aligned}F_i^{MINL}(P_i) &= (a_iP_{i1}^2 + b_iP_{i1} + c_i + |e_i \times \sin(f_i \times (\underline{P_i} - P_{i1}))|)u_{i1} \\ &+ (a_iP_{i2}^2 + b_iP_{i2} + c_i + |e_i \times \sin(f_i \times (\underline{P_i} - P_{i2}))|)u_{i2} \\ &+ (a_iP_{i3}^2 + b_iP_{i3} + c_i + |e_i \times \sin(f_i \times (\underline{P_i} - P_{i3}))|)u_{i3} \\ &= F_{i1}(P_{i1})u_{i1} + F_{i2}(P_{i2})u_{i2} + F_{i3}(P_{i3})u_{i3}\end{aligned} \quad (13)$$

$$u_{i1} + u_{i2} + u_{i3} = 1, \quad u_{ik} \in \{0,1\} \quad (14)$$

where $F_i^{MINL}$ is the mixed-integer nonlinear cost function of unit $i$. In addition, $u_{i1}$, $u_{i2}$, and $u_{i3}$ are the binary decision variables. Moreover, $P_{i1}$, $P_{i2}$, and $P_{i3}$ are the sub-powers, corresponding to sub-functions.

Equations (13) and (14) guarantee the selection of only one sub-function. However, in order to ensure that each sub-function keeps its operating constraints, the sub-power limits (15) and the ramp rate limits (16) are taken into account.

$$\underline{P_{i1}} \le P_{i1} \le \overline{P}_{i1}, \quad \underline{P_{i2}} \le P_{i2} \le \overline{P}_{i2}, \text{ and } \underline{P_{i3}} \le P_{i3} \le \overline{P}_{i3} \quad (15)$$

$$\max(\underline{P_i}, P_i^0 - DR_i) \le P_{i1}u_{i1} + P_{i2}u_{i2} + P_{i3}u_{i3} \le \min(\overline{P}_i, P_i^0 + UR_i) \quad (16)$$

For example, in order to select the third operating zone, the third sub-function must be chosen. Thus, $u_{i3}$ must be set to 1 and based on (14), $u_{i1}$ and $u_{i2}$ are set to 0, which reflects that in (15) and (16), only the limits, related to $P_{i3}$, must be satisfied.

Hence, the general structure of Equations (13)-(16) can be defined as Equations (17)-(20), respectively.

$$F_i^{MINL}(P_i) = \sum_{k=1}^{z_i}(F_{ik}(P_{ik}))u_{ik} \quad (17)$$

$$\sum_{k=1}^{z_i} u_{ik} = 1, \ u_{ik} \in \{0,1\} \tag{18}$$

$$\underline{P_{ik}} \le P_{ik} \le \overline{P}_{ik}, \ k=1,\cdots,z_i \tag{19}$$

$$\max(\underline{P_i},P_i^0 - DR_i) \le \sum_{k=1}^{z_i} P_{ik}u_{ik} \le \min(\overline{P}_i, P_i^0 + UR_i) \tag{20}$$

From Equation (18), it can be deduced that only one of the binary variables ($u_{ik}$) can be adjusted to 1, and the others are adjusted to 0. Therefore, only one of the operating zones and its corresponding sub-function (17) is selected. Consequently, the sub-powers will satisfy their limit (19) and the ramp rate limit (20).

After applying the aforementioned modifications to the ED model, the LB-MINLP framework was obtained as follows:

$$\min \sum_{i=1}^{n_g} F_i^{MINL}(P_i) = \sum_{i=1}^{n_g}\sum_{k=1}^{z_i}((F_{ik}(P_{ik}))\cdot u_{ik}) \tag{21}$$

Subject to:

$$\sum_{i=1}^{n_g}\sum_{k=1}^{z_i} P_{ik}u_{ik} = P_D + P_L \tag{22}$$

$$\underline{P_{ik}} \le P_{ik} \le \overline{P}_{ik}, \ k=1,\cdots,z_i, i=1,\cdots,n_g \tag{23}$$

$$\max(\underline{P_i},P_i^0 - DR_i) \le \sum_{k=1}^{z_i} P_{ik}u_{ik} \le \min(\overline{P}_i, P_i^0 + UR_i), i=1,\cdots,n_g \tag{24}$$

$$\sum_{i=1}^{n_g} S_i \ge S_R \tag{25}$$

$$\sum_{k=1}^{z_i} u_{ik} = 1, \ u_{ik} \in \{0,1\} \tag{26}$$

where the transmission loss (5) and the spinning reserve (8) are modified as (27) and (28), respectively.

$$P_L = \sum_{i=1}^{n_g}\sum_{j=1}^{n_g} \left( (\sum_{k=1}^{z_i} P_{ik}u_{ik}) \cdot B_{ij} \cdot (\sum_{k=1}^{z_j} P_{jk}u_{jk}) \right) + \sum_{i=1}^{n_g} B_{0,i} \cdot \sum_{k=1}^{z_i} P_{ik}u_{ik} + B_{00} \tag{27}$$

$$S_i = \min\{(\overline{P}_i - \sum_{k=1}^{z_i} P_{ik}u_{ik}), \overline{S}_i\}, i=1,\cdots,n_g \tag{28}$$

## 3. Case Studies and Results

The proposed LB-MINLP was tested in two commonly-used test systems (6- and 15-unit test systems) and two large-scale systems (namely the practical 140-unit Korean power system and a 10,000-unit test system).

In order to solve the problems with practical constraints that are highly nonlinear and non-smooth, it is not only essential to create an appropriate model, but also to use a proper solver. In this paper, in order to implement the proposed models, a modeling language for mathematical programming (AMPL) was employed [29]. The nonlinear commercial solver (KNITRO) [30], which is a good solver for mixed-integer nonlinear programming, was used to solve the proposed LB-MINLP model. It should be noted that the simulations were carried out on a 2.67-GHz personal computer, with 3 GB of RAM memory.

### 3.1. Case 1: 6-unit test system

The data of this six-unit test system, containing 26 buses and 46 transmission lines, were obtained from [31]. The load demand was 1,263 MW and for the transmission loss, the B-coefficients were considered with a 100 MVA base capacity. It is worth mentioning that in [31], $B_{00}$ was reported 0.056, while the correct $B_{00}$, used in [31] and other works, was 0.0056 [32]. Each unit of this system had two POZs, resulting in 18 operating zones.

Table 1. Comparison of the Optimal Results of LB-MINLP with Other Approaches for 6-unit System— With POZ, Ramp Rate, and Loss

| Methods | $P_1$ (MW) | $P_2$ (MW) | $P_3$ (MW) | $P_4$ (MW) | $P_5$ (MW) | $P_6$ (MW) | Calc.$P_L$ (MW) | $P_D+P_L$ (MW) | Output (MW) | Viol. (MW) | Cost ($/h) | ACT (p.u.) |
|---|---|---|---|---|---|---|---|---|---|---|---|---|
| TSA [33] | 449.3651 | 182.2520 | 254.2904 | 143.4506 | 161.9682 | 86.0185 | 12.8533 | 1275.8533 | 1277.3448 | -1.4915 | 15,451.63 | 321.9 |
| MSFLA [34] | 449.1444 | 173.0537 | 266.0012 | 127.1123 | 174.2513 | 85.8681 | 13.2571 | 1276.2571 | 1275.4310 | -0.8261 | 15,440.90 | NA |
| DSPSOTSA[33] | 439.2935 | 187.7876 | 261.0260 | 129.4973 | 171.7101 | 86.1648 | 13.1481 | 1276.1481 | 1275.4791 | -0.6690 | 15,441.57 | 20.68 |
| CPSO 2 [35] | 434.4295 | 173.3231 | 274.4735 | 128.0598 | 179.4759 | 85.9281 | 13.3233 | 1276.3233 | 1275.6899 | -0.6334 | 15,446.73 | 245.6 |
| CPSO 1 [35] | 434.4236 | 173.4385 | 274.2247 | 128.0183 | 179.7042 | 85.9082 | 13.3268 | 1276.3268 | 1275.7175 | -0.6093 | 15,447.08 | 241.9 |
| TSA [36] | 451.7300 | 185.2300 | 260.9300 | 133.1000 | 171.0800 | 73.5100 | 13.1519 | 1275.5800 | 1276.1519 | -0.5719 | 15,449.20 | 687.6 |
| MSFL [37] | 445.0140 | 175.5156 | 264.2614 | 137.3012 | 162.7899 | 90.4992 | 12.9434 | 1275.9437 | 1275.3813 | -0.5624 | 15,442.59 | NA |
| AIS [38] | 458.2904 | 168.0518 | 262.5175 | 139.0604 | 178.3936 | 69.3416 | 13.1996 | 1276.1966 | 1275.6553 | -0.5413 | 15,448.00 | NA |
| MPSO [39] | 446.4869 | 168.6612 | 265.0000 | 139.4927 | 164.0036 | 91.7465 | 12.9281 | 1275.9281 | 1275.3909 | -0.5372 | 15,443.09 | NA |
| PSO-SIF [40] | 446.9122 | 173.1470 | 263.6812 | 139.1446 | 165.7765 | 86.7538 | 12.9503 | 1275.9503 | 1275.4155 | -0.5348 | 15,442.66 | 508.6 |
| PC-PSO [41] | 437.7900 | 195.9800 | 256.7200 | 149.3600 | 166.2000 | 69.2600 | 12.8234 | 1275.8234 | 1275.3100 | -0.5134 | 15,453.09 | 2.18 |
| MIPSO [42] | 447.2965 | 173.2582 | 263.6017 | 138.8752 | 165.5300 | 86.8796 | 12.9526 | 1275.9526 | 1275.4412 | -0.5114 | 15,442.98 | NA |
| MIQCQP [24] | 447.4000 | 173.2400 | 263.3800 | 138.9800 | 165.3900 | 87.0500 | 12.9487 | 1275.9487 | 1275.4400 | -0.5087 | 15,443.07 | 13.59 |
| SOH-PSO [41] | 438.2100 | 172.5800 | 257.4200 | 141.0900 | 179.3700 | 86.8800 | 13.0585 | 1276.0585 | 1275.5500 | -0.5085 | 15,446.02 | 2.14 |
| Ө-PSO [43] | 447.1045 | 173.1123 | 263.6503 | 139.1516 | 165.9343 | 86.5037 | 12.9533 | 1275.9533 | 1275.4567 | -0.4966 | 15,443.18 | 26.79 |
| APSO [44] | 446.6686 | 173.1556 | 262.8260 | 143.4686 | 163.9139 | 85.3437 | 12.8617 | 1275.8617 | 1275.3764 | -0.4853 | 15,443.57 | NA |
| PSOM2 [36] | 444.7200 | 172.3700 | 260.5000 | 144.8600 | 167.7100 | 85.2300 | 12.3900 | 1275.8699 | 1275.3900 | -0.4799 | 15,444.50 | 35.88 |
| BA [36] | 438.6500 | 167.9000 | 262.8200 | 136.7700 | 171.7600 | 97.6700 | 13.0495 | 1275.5700 | 1276.0495 | -0.4795 | 15,445.90 | 204.4 |
| IASFLA [45] | 446.7210 | 175.7774 | 264.6118 | 140.2857 | 160.9343 | 87.1002 | 12.8867 | 1275.8867 | 1275.4304 | -0.4563 | 15,442.61 | NA |
| DE [46] | 447.7440 | 173.4070 | 263.4110 | 139.0760 | 165.3640 | 86.9440 | 12.9570 | 1275.9570 | 1275.946 | -0.0110 | 15,449.77 | 0.69 |
| GAAPI [47] | 447.1200 | 173.4100 | 264.1100 | 138.3100 | 166.0200 | 87.0000 | 12.9779 | 1275.9779 | 1275.9700 | -0.0079 | 15,449.81 | NA |
| λ-logic [48] | 447.5076 | 173.3159 | 263.4605 | 139.0629 | 165.4711 | 87.1324 | 12.9580 | 1275.958 | 1275.9504 | - 0.0076 | 15,449.80 | 1.36 |
| TS [49] | 459.0753 | 185.0675 | 264.2094 | 138.1222 | 154.4716 | 74.9900 | 12.9422 | 1275.9422 | 1275.9360 | -0.0062 | 15,454.89 | 624.4 |
| Q-Learning [50] | 448.9480 | 173.5954 | 266.2876 | 127.1212 | 174.3471 | 85.9702 | 13.2740 | 1276.2740 | 1276.2695 | -0.0045 | 15,452.05 | NA |
| MTS [49] | 448.1277 | 172.8082 | 262.5932 | 136.9605 | 168.2031 | 87.3304 | 13.0205 | 1276.0205 | 1276.0231 | +0.0026 | 15,450.06 | 40.52 |
| GA [31] | 474.8066 | 178.6363 | 262.2089 | 134.2826 | 151.9039 | 74.1812 | 13.0217 | 1276.0217 | 1276.0195 | -0.0022 | 15,459.00 | NA |
| PSO [31] | 447.1280 | 173.3221 | 263.4745 | 139.0594 | 165.4761 | 87.1280 | 12.9584 | 1275.9584 | 1275.9571 | -0.0013 | 15,450.00 | NA |
| NPSO-LRS [51] | 446.9600 | 173.3944 | 262.3436 | 139.5120 | 164.7089 | 89.0162 | 12.9361 | 1275.9361 | 1275.9351 | -0.0010 | 15,449.94 | NA |
| PSO-LRS [51] | 447.4440 | 173.3430 | 263.3646 | 139.1279 | 165.5076 | 87.1698 | 12.9571 | 1275.9571 | 1275.9569 | -0.0002 | 15,449.90 | NA |
| NPSO [51] | 447.4734 | 173.1012 | 262.6804 | 139.4156 | 165.3002 | 87.9761 | 12.9470 | 1275.9470 | 1275.9469 | -0.0001 | 15,449.91 | NA |
| PSO [49] | 447.5823 | 172.8387 | 261.3300 | 138.6812 | 169.6781 | 85.8963 | 13.0066 | 1276.0066 | 1276.0066 | 0.0000 | 15,450.14 | 207.2 |
| LB-MINLP | 447.5038 | 173.3182 | 263.4628 | 139.0653 | 165.4734 | 87.1347 | 12.9582 | 1275.9582 | 1275.9582 | 0.0000 | 15,449.89 | 1.00 |

Table 1 presents the comparison of the proposed LB-MINLP model and the other thirty-one algorithms, including artificial immune system (AIS) [38], adaptive particle swarm optimization (APSO) [44], bees algorithm (BA) [36], particle swarm optimization (PSO) [31], [49], chaotic PSO (CPSO) [40], differential evolution (DE) [46], tabu search (TS) [49], tabu search algorithm (TSA)

[36], [33], distributed Sobol PSO and TSA (DSPSO-TSA) [33], genetic algorithm (GA) [31], genetic algorithm and ant colony optimization (GAAPI) [47], improved adaptive shuffled frog-leaping algorithm [45], modified iterative PSO (MIPSO) [42], mixed-integer quadratically constrained quadratic programming (MIQCQP) [24], modified PSO (MPSO) [39], modified shuffled frog-leaping (MSFL) [37], modified shuffled frog-leaping algorithm (MSFLA) [34], multiple tabu search (MTS) [49], new PSO (NPSO) [51], PSO with local random search (PSO-LRS) [51], NPSO with local random search (NPSO-LRS) [51], passive congregation-based PSO (PC-PSO) [41], PSO with mutation operators (PSOM) [36], PSO with smart inertia factor (PSO-SIF) [40], Q-learning [50], self-organizing hierarchical PSO (SOH-PSO) [41], efficient lambda logic based optimization procedure (λ-logic) [48], and Θ-PSO [43]. In this table, the transmission loss ($P_L$) was re-calculated based on the reported outputs of the generators in each work, while the total demand and the loss ($P_D + P_L$) are reported to check the power balance and the probable violations. It is worth stating that in order to compare the CPU time of different approaches, the adjusted CPU time was adopted [52].

$$ACT = \frac{\text{given CPU speed (GHz)}}{2.67 \text{ GHz}} \times \frac{\text{given CPU time (s)}}{\text{CPU time from LB-MINLP (s)}} \tag{29}$$

In Table 1, the algorithms were sorted according to their violation in descending order, so that the algorithms with the highest violation are at the top of the table. From this table, it can be realized that the power balance has been violated in 96.77% (30 out of 31) of the works in the literature. Hence, in the current highly competitive world, such a solution is defective. By comparing the results of the proposed LB-MINLP model with the PSO algorithm expressed in [49], which is the only unviolated approach in the literature, it can be seen that the LB-MINLP model has obtained a better solution ($15,449.89/h). In addition, the superiority of the LB-MINLP model is not only in finding the best solution with a perfect precision (without violation), but also its computational efficiency is remarkably high (only 0.031 seconds), which is nearly 207 times faster than the only unviolated approach among all, PSO, [49]. Therefore, the proposed MINLP-based model acts better than all the other works in the literature and this finding proves the suggested claim that the proposed nonlinear non-convex model is a solver-friendly model and does not bring difficulties for the commercial solver. Moreover, as mentioned in the introduction section, the MIQCQP approach, which is the only solver-based model, shows a considerable violation, where the total generated power is 0.5087 less than the total system demand (it means, the consumer demand plus the transmission loss).

### 3.2. Case 2: 15-unit test system

The 15-unit test system is one of the most commonly used systems in the literature that has been so far tested under different conditions and topologies. Therefore, this case was employed to verify the proposed model and check its velocity in comparison with the other approaches in the literature. The generating units 2, 5, 6, and 12, among the 15 units, had a total of 11 prohibited zones.

### 3.2.1. First Condition: With prohibited operating zones, ramp rate, and without losses

The 15-unit test system, by considering the operating zones and without considering the transmission losses, has been widely studied to validate the existing approaches in the literature and mainly the models that are solvable via commercial solvers. We used this case to show the accuracy and the computational efficiency of the proposed mixed-integer nonlinear model, in comparison with the heuristic-based and linear models. The load demand was assigned to 2,650 MW, with the system spinning reserve requirement of 200 MW. The data of this system were obtained from [53].

Table 2. Comparison the Optimal Results of LB-MINLP with Other Approaches for 15-unit System— With Ramp Rate and without Loss

| Method | Cost ($/h) | Time (s) | ACT (p.u.) |
|---|---|---|---|
| QEA [54] | 32,507.48 | NA* | NA |
| IFEP [55] | 32,507.46 | 3.138 | NA |
| ESO [56] | 32,506.6 | 13.79 | NA |
| PSO-TVAC [57] | 32,506.43 | NA | NA |
| PSO [58] | 32,506.3 | 1.969 | NA |
| $\lambda$-logic [48] | 32,506.18 | 0.032 | 2.61 |
| IQEA [54] | 32506.14 | NA | NA |
| MIQP [21] | 32506.14 | 0.25 | 17.024 |
| M-MIQP [22] | 32506.14 | 0.03 | 2.043 |
| UDB-MIQP [23] | 32506.14 | 0.007 | 0.636 |
| LB-MINLP | 32506.14 | 0.011 | 1.000 |

Table 2 exposes the comparison of the proposed LB-MINLP and the other approaches in the literature. Among the heuristic-based methods, the quantum-inspired evolutionary algorithm (QEA) [54], improved fast evolutionary programming (IFEP) [55], efficient evolutionary strategy optimization (ESO) [56], PSO with time-varying acceleration coefficients (PSO-TVAC) [57], PSO [58], and $\lambda$-logic [48] obtained near-optimal solutions, while the improved QEA (IQEA) [54] achieved the global optimal solution, as reported by the MIQP-based model in [5, 6, 44]. However, the results acquired by the proposed LB-MINLP model not only exhibited the required accuracy of the model for finding the global optimal solution, but also confirmed the high computational efficiency compared with the quadratic-based models, which are well-known for high computational efficiency. Our suggested nonlinear model also displayed a much higher computational efficiency over the MIQP and the M-MIQP, by approximately 17.024 and 2.043 times higher velocity, respectively. Having compared the UDB-MIQP, which is a well-defined user-friendly MIQP-based model, the proposed LB-MINLP model required more CPU time (about

0.004 s). However, this negligible increment in the CPU time is the price that the LB-MINLP model must pay for the non-convexity and nonlinearity of addressing the shortcomings of the MIQP-based models and for considering the valve-point effects and the transmission losses. These findings demonstrated the effectiveness of the LB-MINLP model in finding the best optimal solution, while the appropriate trade-off between the complexity and the accuracy of the model resulted in a high computational efficiency.

### 3.2.2. *Second Condition: With prohibited operating zones, ramp rate, and transmission losses*

Under this condition, in addition to considering the POZs, the transmission loss was also taken into account. The load demand was assigned to 2,630 MW, with a system spinning reserve requirement of 200 MW. The generators' data and the B-coefficients of the 15-unit test system were gained from [31] with respect to the change reported in [32]. It is worth mentioning that the costs and losses presented in the tables were re-calculated based on the outputs of the reported units, and the presented approaches were sorted in descending order, according to their violation.

Table 3. Comparison the Generators' Outputs of LB-MINLP with Other Approaches for 15-unit System—With Ramp Rate and Loss

| Output power | MIQCQP [24] | CCPSO [59] | BA [60] | FA [61] | LB-MINLP |
|---|---|---|---|---|---|
| $P_1$ (MW) | 455.0000 | 455.0000 | 455.0000 | 455.0000 | 455.0000 |
| $P_2$ (MW) | 380.0000 | 380.0000 | 380.0000 | 380.0000 | 380.0000 |
| $P_3$ (MW) | 130.0000 | 130.0000 | 130.0000 | 130.0000 | 130.0000 |
| $P_4$ (MW) | 130.0000 | 130.0000 | 130.0000 | 130.0000 | 130.0000 |
| $P_5$ (MW) | 170.0000 | 170.0000 | 170.0000 | 170.0000 | 170.0000 |
| $P_6$ (MW) | 460.0000 | 460.0000 | 460.0000 | 460.0000 | 460.0000 |
| $P_7$ (MW) | 430.0000 | 430.0000 | 430.0000 | 430.0000 | 430.0000 |
| $P_8$ (MW) | 72.1300 | 71.7526 | 71.7474 | 71.7450 | 71.7430 |
| $P_9$ (MW) | 58.5400 | 58.9090 | 58.9140 | 58.9164 | 58.9184 |
| $P_{10}$ (MW) | 160.0000 | 160.0000 | 160.0000 | 160.0000 | 160.0000 |
| $P_{11}$ (MW) | 80.0000 | 80.0000 | 80.0000 | 80.0000 | 80.0000 |
| $P_{12}$ (MW) | 80.0000 | 80.0000 | 80.0000 | 80.0000 | 80.0000 |
| $P_{13}$ (MW) | 25.0000 | 25.0000 | 25.0000 | 25.0000 | 25.0000 |
| $P_{14}$ (MW) | 15.0000 | 15.0000 | 15.0000 | 15.0000 | 15.0000 |
| $P_{15}$ (MW) | 15.0000 | 15.0000 | 15.0000 | 15.0000 | 15.0000 |
| $P_L$ (MW) | 30.6635* | 30.6615* | 30.6614 | 30.6614 | 30.6614 |
| $P_D+P_L$ (MW) | 2660.67* | 2660.6615* | 2660.6614 | 2660.6614 | 2660.6614 |
| Cost ($/h) | 32,704.53* | 32,704.45 | 32,704.45 | 32,704.45 | 32,704.45 |

* The results with modification based the reported units' outputs.

Table 3 discloses the optimal output of all units of the 15-unit test system, by considering the ramp rate and the transmission losses, and provides its comparison with four different approaches, which obtained the best optimal solution reported. These approaches are MIQCQP [24], PSO with chaotic sequences and crossover operation (CCPSO) [59], firefly algorithm (FA) [61], and bat algorithm (BA) [60]. The outputs of the units showed that there only exists a little difference between these approaches, in the outputs of units 8 and 9, which can be the consequence of the model accuracy or the solution methodology. The MIQCQP and CCPSO methods uncovered a small violation in calculating the losses and the costs. However, in order to make it clear and show the effectiveness of the proposed LB-MINLP in facing the transmission losses, the optimal solutions of different approaches were considered in more depth in Table 4, including chaotic PSO (CPSO) [35], hybrid swarm intelligence based harmony search algorithm (HHS) [62], artificial bee colony algorithm (ABC) [63], modified DE (MDE) [64], bat algorithm (BA) [60], and the other approaches, named in the previous sections.

Table 4. Comparison of the Proposed LB-MINLP with Other Approaches for 15-unit System— With POZ, Ramp Rate, and Loss

| Methods | Reported $P_L$ (MW) | Calculated $P_L$ (MW) | $P_D+P_L$ (MW) | Output (MW) | Violation (MW) | Cost ($/h) | Time (s) | ACT (p.u.) |
|---|---|---|---|---|---|---|---|---|
| CPSO 1 [35] | 32.1302 | 36.1371 | 2666.1371 | 2662.0520 | -4.0851 | 32,835.00 | 13.31 | 90.31 |
| CPSO 2 [35] | 32.1303 | 36.1396 | 2666.1396 | 2662.1000 | -4.0396 | 32,834.00 | 13.13 | 89.09 |
| MPSO [39] | 29.9780 | 29.8846 | 2659.8846 | 2661.6235 | +1.7389 | 32,738.42 | NA | NA |
| APSO [44] | 28.3700 | 29.5866 | 2659.5866 | 2658.3100 | -1.2766 | 32,724.78 | NA | NA |
| HHS [62] | 29.6631 | 30.6559 | 2660.6559 | 2659.6594 | -0.9965 | 32,692.86 | 4.7231 | 25.64 |
| GAAPI [47] | 30.3600 | 29.6332 | 2659.6332 | 2660.3700 | +0.7368 | 32,732.95 | NA | NA |
| MIPSO [42] | 30.1000 | 30.6397 | 2660.6397 | 2660.0556 | -0.5841 | 32,697.54 | NA | NA |
| MTS [49] | 31.3523 | 31.0647 | 2661.0647 | 2661.3635 | +0.2988 | 32,716.87 | 3.65 | 25.76 |
| ABC [63] | 30.9591 | 30.6775 | 2660.6775 | 2660.9591 | +0.2816 | 32,707.85 | 11.02 | 50.84 |
| IASFLA [45] | 29.9374 | 30.7534 | 2660.7534 | 2660.9609 | +0.2075 | 32,712.03 | NA | NA |
| PSO-SIF [40] | 30.8822 | 30.6848 | 2660.6848 | 2660.8822 | +0.1974 | 32,706.88 | 60.7 | 329.5 |
| IPSO [65] | 30.8574 | 30.6739 | 2660.6739 | 2660.8574 | +0.1835 | 32,706.66 | 2.36 | 19.60 |
| SPSO [41] | 30.4900 | 30.2605 | 2660.2605 | 2660.4400 | +0.1795 | 32,798.69 | 0.0913 | 0.69 |
| Ө-PSO [43] | 30.8260 | 30.6504 | 2660.6504 | 2660.8270 | +0.1766 | 32,706.55 | 5.5794 | 37.86 |
| MSFL [37] | 30.8570 | 30.6733 | 2660.6733 | 2660.8497 | +0.1764 | 32,706.57 | NA | NA |
| PSO [32] | 37.3329 | 37.3329 | 2667.3329 | 2667.4250 | +0.0921 | 33,020.00 | 26.59 | NA |
| TSA [33] | 33.8110 | 33.8104 | 2663.8104 | 2663.7200 | -0.0904 | 32,918.00 | 25.75 | 111.8 |
| MDE [64] | 30.6200 | 30.6169 | 2660.6169 | 2660.5500 | -0.0669 | 32,704.90 | 12.88 | 97.88 |
| AIS [38] | 32.4075 | 32.3452 | 2662.3452 | 2662.4040 | +0.0588 | 32,854.00 | NA | NA |
| SOH-PSO [41] | 32.2800 | 32.2340 | 2662.2340 | 2662.2900 | +0.0560 | 32,751.39 | 0.0936 | 0.71 |
| PC-PSO [41] | 30.5400 | 30.5915 | 2660.5915 | 2660.5500 | -0.0415 | 32,775.36 | 0.0967 | 0.73 |
| GA [32] | 38.2499 | 38.2499 | 2668.2499 | 2668.2900 | +0.0401 | 33,149.00 | 49.31 | NA |
| DSPSO–TSA [33] | 30.9520 | 30.9520 | 2660.9520 | 2660.9620 | +0.0100 | 32,715.06 | 2.30 | 9.99 |
| MIQCQP [24] | 30.6600 | 30.6635 | 2660.6634 | 2660.6700 | -0.0065 | 32,704.53 | 4.65 | 31.55 |
| CCPSO [59] | 30.6616 | 30.6615 | 2660.6115 | 2660.6116 | +0.0001 | 32,704.45 | 16.20 | 87.93 |
| λ-logic [48] | 29.9491 | 29.9491 | 2659.9491 | 2659.9491 | 0.0000 | 32,713.95 | 0.016 | 0.11 |
| BA [60] | 30.6614 | 30.6614 | 2660.6114 | 2660.6614 | 0.0000 | 32,704.45 | NA | NA |
| FA [61] | 30.6614 | 30.6614 | 2660.6114 | 2660.6614 | 0.0000 | 32,704.45 | 16.05 | 95.83 |
| LB-MINLP | 30.6614 | 30.6614 | 2660.6614 | 2660.6614 | 0.0000 | 32,704.45 | 0.138 | 1.00 |

Table 4 exhibits that the λ-logic [48], BA [60], FA [61], and the proposed LB-MINLP have found near-optimal/optimal solutions without violation in the power balance. Added to them, the MIQCQP model, which is a solver-based model, the same as the 6-unit case, presented a violation in the optimal solution, with a probable failure in the bi-level procedure or a miscalculation of the transmission loss. It was also comprehended that the LB-MINLP model has obtained the best optimal solution reported in the literature (over the BA and the FA approaches), while it is much faster (96 times) than the FA approach. However, the λ-logic failed to make an appropriate trade-off between the preciseness and the computational efficiency by obtaining a low-quality solution within a perfect CPU time.

### 3.2.3. *Third Condition: With prohibited operating zones, ramp rate, valve-point effect, and losses*

The valve-point effect was taken into account for all generators in order to test the LB-MINLP model under a more complex condition. The data of this system was the same as the previous case, and the valve-point data can be observed in the Appendix (Table A1). The main reasons for presenting this highly nonlinear and non-smooth test system are as follows: 1) To show that even though the simultaneous consideration of the POZs and the valve-point effects in generation units (in this case: units 2, 5, 6, and 12) brings some difficulties in solving the process and results in increasing the CPU time, the proposed LB-MINLP model ultimately conquers it and finds an optimal solution; 2) To provide a critical system for future researchers to investigate the efficiency of their approaches.

Table 5. Optimal Result of LB-MINLP for 15-unit System— With Loss, Ramp Rate, Reserve, and Valve Point Effect

| Unit | Power (MW) | Unit | Power (MW) | Unit | Power (MW) | $P_L$ (MW) | Cost ($/h) | Time (s) |
|---|---|---|---|---|---|---|---|---|
| 1 | 455.00 | 6 | 456.53 | 11 | 58.31 | | | |
| 2 | 380.00 | 7 | 417.38 | 12 | 55.00 | | | |
| 3 | 100.55 | 8 | 141.60 | 13 | 25.00 | 38.71 | 33,417.08 | 12.425 |
| 4 | 100.55 | 9 | 105.56 | 14 | 55.00 | | | |
| 5 | 150.00 | 10 | 153.23 | 15 | 15.00 | | | |

It can be witnessed in Table 5 that unlike the previous test, in which the optimal solution was achieved in only 0.138 s, the solution of this case was obtained in 12.425 s. This proves the difficulties of finding an optimal solution for such critical cases. However, unlike the other two conditions, where all the units (except units 8 and 9) were generating at their maximum or minimum limits, under this condition, most of the units were generating under the limit, which showed the difficulties of a solver in dispatch adjustments.

Table 6. Operating zones of the generating units under different conditions via LB-MINLP for 15-unit System

| unit | First Condition | | | | Second Condition | | | | Third Condition | | | |
|---|---|---|---|---|---|---|---|---|---|---|---|---|
| | Zone 1 | Zone 2 | Zone 3 | Zone 4 | Zone 1 | Zone 2 | Zone 3 | Zone 4 | Zone 1 | Zone 2 | Zone 3 | Zone 4 |
| 2 | | | | ● | ● | | | ● | | | ● | |
| 5 | | ● | | | ● | | | | ● | | | |
| 6 | | | | ● | | | | ● | | | | ● |
| 12 | | | ● | | | | | ● | | | ● | |

Table 6 displays the operating zone, in which the units with POZs are operating. As can be seen, by comparing with the first condition, the obtained pattern indicated only one change for the second condition (only in unit 5), while under the third condition, the patterns of three out of four units (units 1, 2, and 12) changed. This confirmed the difficulties that a simultaneous consideration of POZ, ramp rate, valve point, and transmission loss may bring for the solvers.

### 3.3. Case 3: Large-scale power systems: 140- and 10,000-unit test systems

In order to demonstrate the effectiveness of the proposed LB-MINLP model in solving large-scale power systems, a practical Korean 140-unit power system and a 10,000-unit test system were employed. The data of the Korean 140-unit were gained from [59], and the data of the 10,000 unit, containing two hundred and fifty 40-unit systems, with 11,500 POZs, were the same as the data stated in [21], and the load demand was assigned to 1,750,000 MW.

Table 7. Comparison of the Optimal Results of LB-MINLP with Other Works for 140-unit System— With POZ, Ramp Rate, and Valve Point

| Method | Cost ($/h) | Standard Deviation | Time (s) | Violation (MW) | ACT (p.u.) |
|---|---|---|---|---|---|
| MFPA [66] | 1,657,962.69 | 0.05557 | 5.71 | 0.0001 | NA |
| DEL [67] | 1,657,962.72 | 57.9800 | NA | -----* | NA |
| CTPSO [59] | 1,657,962.73 | 7.3150 | 100.00 | NA | 24.0 |
| CCEDE [68] | 1,657,962.73 | 1.1466 | NA | 0.0000 | NA |
| CQGSO [69] | 1,657,962.73 | 0.0661 | 31.67 | NA | NA |
| CSPSO [59] | 1,657,962.73 | 0.0235 | 99.00 | NA | 23.8 |
| COPSO [59] | 1,657,962.73 | 0.0002 | 150.00 | NA | 36.0 |
| DSD [20] | 1,657,962.73 | NA | 0.177** | NA | 0.1 |
| CCPSO [59] | 1,657,962.73 | 0.0000 | 150.00 | NA | 36.0 |
| LB-MINLP | 1,657,962.72 | 0.0000 | 2.12 | 0.0000 | 1.0 |

* There is a mistyping in presenting of the unit's output.
** There is no information about the time spent to adjust the model coefficients that directly affects the solution.

Table 7 presents the comparison of the results obtained from the LB-MINLP model and the outcome of the other approaches in the literature. The gathered results exhibited that the LB-MINLP model can obtain the optimal solution (without violation and standard deviation) in only 2.12 s, which is fast enough for this large-scale practical system. Among the presented approaches, the CCPSO is the only approach without standard deviation, while for the DSD model, there is no information about the standard deviation. The performance of the DSD approach, which is the quickest approach among all, highly depends on the coefficient adjustment and according to the authors' claim, adjusting the coefficient of the DSD approach is a trial-and-error process, which is a significant drawback in real-world applications. On the other hand, in the DSD approach, the simplification of the transmission loss formula decreases the model accuracy and this would be the price of increasing the velocity. However, the LB-MINLP model addresses such drawbacks by introducing a solver-friendly framework, which is solvable by commercial solvers and can be easily applied to different systems, regardless the system topology or the operating conditions.

The proposed LB-MINLP model was tested on a very large-scale 10,000-unit power system to determine its performance. The optimal solution of this system was $25,191,921.7998/h, which is five times of the 2,000-unit system, expressed in [23]. This proved the accuracy and the

effectiveness of the LB-MINLP model in solving a very large-scale system (with thousands of POZs) in an acceptable CPU time of 23.39 seconds.

## 4. Conclusion

In this paper, a logic-based mixed-integer nonlinear programming framework (LB-MINLP) was proposed to solve non-convex and non-smooth economic dispatch (ED) problems. The LB-MINLP model addressed the existing gap of having an accurate solver-based model for solving the ED problems by considering the non-smooth and non-convex terms, transmission losses, and prohibited operating zones. The proposed solver-friendly model perfectly addressed the current drawbacks in the literature (for heuristic- or solver-based models), for considering the losses in non-smooth ED problems with disjoint operating zones. The acquired results confirmed the superiority of the proposed LB-MINLP model in solving the ED problems so that it is not only capable of finding the most precise solution, but also is fast enough to be applied in online-based problems. Moreover, the results verified the effectiveness and usefulness of the LB-MINLP model in finding the optimal solution of very large-scale power systems.

Some prospects for the future works can be to investigate an equivalent and modified model for solving the operation- and the planning-based problems, by considering disjoint operating zones (for example, optimal power flow and unit commitment).

## 5. Acknowledgements

This work was supported by FAPESP (Nos. 2014/22828-3, and 2016/14319-7), CAPES, and CNPq (No. 305371/2012-6).

## 6. Appendix

In this appendix, the valve point coefficients of the 15-unit test system are presented within Table A1.

Table A1. Valve Point Coefficients of the 15-unit Test System

| Unit | $e$ | $f$ | Unit | $e$ | $f$ | Unit | $e$ | $f$ |
|---|---|---|---|---|---|---|---|---|
| 1 | 170 | 0.091 | 6 | 120 | 0.035 | 11 | 110 | 0.082 |
| 2 | 110 | 0.078 | 7 | 100 | 0.089 | 12 | 450 | 0.089 |
| 3 | 275 | 0.039 | 8 | 230 | 0.077 | 13 | 200 | 0.042 |
| 4 | 275 | 0.039 | 9 | 250 | 0.039 | 14 | 150 | 0.063 |
| 5 | 120 | 0.077 | 10 | 200 | 0.049 | 15 | 175 | 0.045 |